\documentclass[11pt]{lmcs}
\pdfoutput=1

\usepackage{lastpage}
\lmcsdoi{16}{2}{2}
\lmcsheading{}{\pageref{LastPage}}{}{}%
{Apr.~16,~2018}{Apr.~24,~2020}{}

\usepackage{amsmath,amssymb,amscd,amsthm,amsfonts,amstext,amsbsy,mathrsfs,upgreek,mathtools,stmaryrd,enumitem,bbm}

\newcommand{\NN}{\mathbb{N}}

\newcommand{\cof}{\operatorname{cof}}

\newcommand{\Card}{{\mathrm{Card}}}

\newcommand{\ran}[1]{{{\rm{ran}}(#1)}}

\newcommand{\ITTM}{$\mathrm{ITTM}$}
\newcommand{\UU}{U}
\newcommand{\Ord}{\mathrm{Ord}}

\newenvironment{enumerate-(a)}{\begin{enumerate}[label={\upshape (\alph*)}]}{\end{enumerate}}

\newenvironment{enumerate-(a)-r}{\begin{enumerate}[label={\upshape (\alph*)},resume]}{\end{enumerate}}

\newenvironment{enumerate-(A)}{\begin{enumerate}[label={\upshape (\Alph*)}]}{\end{enumerate}}

\newenvironment{enumerate-(A)-r}{\begin{enumerate}[label={\upshape (\Alph*)},resume]}{\end{enumerate}}

\newenvironment{enumerate-(i)}{\begin{enumerate}[label={\upshape (\roman*)}]}{\end{enumerate}}

\newenvironment{enumerate-(i)-r}{\begin{enumerate}[label={\upshape (\roman*)},resume]}{\end{enumerate}}

\newenvironment{enumerate-(I)}{\begin{enumerate}[label={\upshape (\Roman*)}]}{\end{enumerate}}

\newenvironment{enumerate-(I)-r}{\begin{enumerate}[label={\upshape (\Roman*)},resume]}{\end{enumerate}}

\newenvironment{enumerate-(1)}{\begin{enumerate}[label={\upshape (\arabic*)}]}{\end{enumerate}}

\newenvironment{enumerate-(1)-r}{\begin{enumerate}[label={\upshape (\arabic*)},resume]}{\end{enumerate}}

\begin{document}




\author{Merlin Carl} 
\address{Fachbereich Mathematik und Statistik,
  University of Konstanz,
  78457 Konstanz,
  Germany, and
  Europa-Universit{\"a}t Flensburg,
  Institut f{\"u}r mathematische,
  naturwissenschaftliche und technische Bildung,
  Abteilung f{\"u}r Mathematik und ihre Didaktik,
  Auf dem Campus 1b, 24943 Flensburg, Germany}
\email{merlin.carl@uni-konstanz.de} 

\author{Benjamin Rin}
\address{Departement Filosofie en Religiewetenschap, 
Utrecht University, 
Janskerkhof 13, 
3512 BL, Utrecht, 
The Netherlands}
\email{b.g.rin@uu.nl} 

\author{Philipp Schlicht}
\address{Department of Computer Science, 
The University of Auckland, 
Private Bag 92019, Auckland 1142, 
New Zealand, and 
School of Mathematics, 
University of Bristol, 
Fry Building.
Woodland Road, 
Bristol, BS8~1UG, UK}
\email{philipp.schlicht@bristol.ac.uk}

\thanks{We would like to thank the anonymous referees for their helpful comments. 
This project has received funding from the European Union's Horizon 2020 research and innovation programme under the Marie Sk\l odowska-Curie grant agreement No 794020 (IMIC) for the third-listed author.
} 

\date{\today}

\title[Reachability for infinite time Turing machines with long tapes]{Reachability for infinite time Turing machines with long tapes}


\begin{abstract} 
Infinite time Turing machine models with tape length $\alpha$, denoted $T_\alpha$, 
strengthen the machines of Hamkins and Kidder 
with tape length $\omega$. 
A new phenomenon is that for some countable ordinals $\alpha$, some cells cannot be halting positions of $T_\alpha$ given trivial input. 
The main open question in 
a paper of Rin from 2014 
asks about the size of the least such ordinal $\delta$. 

We answer this by providing various characterizations. 
For instance, $\delta$ is the least ordinal with any of the following properties: 
\begin{itemize} 
\item 
For some $\xi<\alpha$, there is a $T_\xi$-writable but not $T_\alpha$-writable subset of $\omega$. 
\item 
There is a gap in the $T_\alpha$-writable ordinals. 
\item 
$\alpha$ is uncountable in $L_{\lambda_\alpha}$. 
\end{itemize} 
Here $\lambda_\alpha$ denotes the supremum of $T_\alpha$-writable ordinals, i.e. those with a $T_\alpha$-writable code of length $\alpha$. 

We further use the above characterizations, and an analogue to Welch's submodel characterization of the ordinals $\lambda$, $\zeta$ and $\Sigma$, to show that $\delta$ is large in the sense that it is a closure point of the function $\alpha \mapsto \Sigma_\alpha$, where $\Sigma_\alpha$ denotes the supremum of the $T_\alpha$-accidentally writable ordinals. 
\end{abstract} 

\maketitle

\setcounter{tocdepth}{2}

\section{Introduction}

\subsection{Motivation} 

The infinite time Turing machines  introduced by Hamkins and Kidder (see \cite{HL}) are, roughly, Turing machines with a standard tape that run for transfinite ordinal time. 
One of the main motivations for studying these machines is the fact that 
they model a class of functions that is closely related to $\Sigma^1_1$ and $\Pi^1_1$ sets in descriptive set theory. 
Soon after they were introduced, several variations were proposed, for instance with an arbitrary ordinal as tape length \cite{K}, 
or an exponentially closed ordinal as tape length and time bound \cite{K1,KS}. 

More recently, the second author studied machines with an arbitrary ordinal as tape length, but no ordinal bound on the running time \cite{MR3210081}. 
These machines are natural generalizations of infinite time Turing machines for tapes of length $\alpha$, and are thus called \emph{$\alpha$-\ITTM s}.
They do not include ordinal parameters, which are present in most other models \cite{K1,KS,CO}.
For a given ordinal $\alpha \in \text{On}$ and program $e \in \mathbb N$ there exists a unique machine, $T_\alpha [e]$. 
We will frequently identify a program with the corresponding machine. 
The set $\{T_\alpha [e] \mid e \in \mathbb N\}$ of all machines with tape length $\alpha$ is called the \emph{device} or \emph{machine model} $T_\alpha$. It is 
known that the computability 
strength\footnote{Computability strength in the present sense is a relative notion: given $\alpha,\beta \in \text{On}$, we write $T_\alpha \preceq T_\beta$ when the set of functions  $f\colon 2^{\min(\alpha,\beta)} \to 2^{\min(\alpha,\beta)}$ computable by $T_\alpha$ (that is to say, computable by $T_\alpha [e ]$ for some $e \in \mathbb N$) is a subset of the set of such functions computable by $T_\beta$.} 
of $T_\alpha$ can increase with $\alpha$, though it remains equal when the increase in $\alpha$ is small. When $\alpha$ itself is not too large, increasing its size necessarily makes the computational strength greater or equal. However, it turns out that for sufficiently large tapes, the machine models' computability strengths are not always commensurable: 
there exist pairs of countable ordinals such that two devices 
with these tape lengths can each compute functions that the other one can't \cite[Proposition 2.9]{MR3210081}. 
Thus $\alpha$-ITTMs fail to be linearly ordered by computational strength. What is responsible for this phenomenon is that, in spite of the lack of ordinal parameters, a machine can use its tape length $\alpha$ to perform computations that rely on the exact size of $\alpha$---an ability which, \emph{because} of the lack of parameters, can permit two differently sized machines to exploit their tape lengths in ways the other cannot. This phenomenon clearly does not occur for models that include ordinal parameters (as in \cite{CO}), since one can then always simulate a shorter tape on a longer one (cf.~\cite[Proposition 2.1]{MR3210081}). This is because one can easily move the head up to cell $\xi$ and halt there whenever one is allowed to mark the $\xi$th cell, as is possible when computing with parameter $\xi$. Indeed, it is straightforward to see that an ordinal parameter $\xi < \alpha$ is equivalent to an oracle that allows a machine $T_\alpha$ to emulate the computational behavior of smaller machine $T_\xi$.

In the present article, we are interested in the writability strength of $\alpha$-\ITTM s without parameters, i.e., in the set of possible outputs of such a machine at the time when it halts. 
One of the tools in \cite{MR3210081} to help classify the machines in question 
is the connection between computability strength and ordinals $\alpha$ such that $T_\alpha$ 
cannot reach\footnote{An ordinal $\mu$ is defined to be \emph{reachable} by $T_\alpha$ when 
there exists a program $P$ such that $T_\alpha$ running $P$ on trivial input (input $\vec 0$) halts with the final 
head position located at cell $\mu$.} all of its cells. 
In particular, let $\delta$ denote the least such ordinal. 
Rin already showed that $\delta$ equals the least ordinal $\gamma$ such that the computability strength of $T_\gamma$ (as defined above) is incomparable with that of some machine with a shorter 
tape \cite[Proposition 2.1]{MR3210081}. 
The main question left open was to identify $\delta$.\footnote{See the discussion after \cite[Proposition 2.9]{MR3210081}.}  
We answer this by giving various characterizations of $\delta$ in the next theorem. 
Some of them are formulated via $\alpha$-\ITTM s. The remaining ones are stated in term of constructible set theory and resemble fine-structural properties of the constructible universe $L$, where first-order definability is replaced with variants of infinite time writability. 

\begin{thm} \label{main theorem} 
The next conditions for $\alpha>\omega$ hold for the first time at the same ordinal: 
\begin{enumerate}[(a)]
\item \label{main theorem 1}
Not every cell is $T_\alpha$-reachable ($T_\alpha$-eventually reachable). 
\item \label{main theorem 2}
There is a gap in the $T_\alpha$-writable ($T_\alpha$-eventually writable) ordinals. 
\item \label{main theorem 3} 
For some $\mu,\nu$ with $\omega\leq\mu\leq\nu<\alpha$, there is a $T_\nu$-writable but not $T_\alpha$-writable subset of $\mu$. 
\item \label{main theorem 4}
$\lambda_\alpha<\hat{\lambda}_\alpha$. 
\item \label{main theorem 5}
$\zeta_\alpha<\hat{\zeta}_\alpha$. 
\item \label{main theorem 6}
$\alpha$ is uncountable in $L_{\lambda_\alpha}$. 
\item \label{main theorem 7}
$\alpha$ is regular in $L_{\lambda_\alpha}$. 
\item \label{main theorem 8}
$\alpha$ is a cardinal in $L_{\lambda_\alpha}$. 
\item \label{main theorem 9}
As \ref{main theorem 6}, \ref{main theorem 7} or \ref{main theorem 8}, but for $\hat{\lambda}_\alpha$, $\zeta_\alpha$, $\hat{\zeta}_\alpha$ or $\hat{\Sigma}_\alpha=\Sigma_\alpha$. 
\end{enumerate} 
\end{thm} 
\noindent
In \ref{main theorem 4}--\ref{main theorem 9}, $\lambda_\alpha$ and $\zeta_\alpha$ denote the versions of $\lambda$ and $\zeta$ for $T_\alpha$ without ordinal parameters, while $\hat{\lambda}_\alpha$ and $\hat{\zeta}_\alpha$ denote those with parameters. 

Let $\lambda$, $\zeta$ and $\Sigma$ denote the suprema of writable, eventually writable and accidentally writable ordinals for \ITTM s. 
The previous characterizations imply that $\Sigma<\delta$. 
By \ref{main theorem 6}, we can obtain triples $(\mu,\nu,\xi)$ with $\mu<\nu<\xi<\delta$ and $L_\mu\prec L_\nu\prec L_\xi$ by forming countable elementary substructures of $L_\delta$ in $L_{\hat{\lambda}_\delta}$. 
Thus $\Sigma<\delta$ holds by Welch's submodel characterization of $\lambda$, $\zeta$ and $\Sigma$ (cf.~\cite[Theorem 30 \& Corollary 32]{MR2493990}). 

The next result (cf. Theorem \ref{Sigmaalpha is less than delta}) is proved 
via a variant of the submodel characterization for $\alpha$-\ITTM s (cf. Theorem \ref{lambda-zeta-Sigma submodel characterization}). 

\begin{thm} \label{bound on delta} 
$\Sigma_{\xi}<\delta$ for all $\xi<\delta$.\footnote{This strengthens the result from \cite{MR3210081} that $\zeta<\delta$ and an unpublished result by Robert Lubarsky that $\Sigma<\delta$. } 
\end{thm}

The structure of the paper is as follows. 
Section \ref{section setting} contains some background on $\alpha$-\ITTM s. 
In Section \ref{section writable and clockable ordinals}, we prove some auxiliary results about writable 
and clockable ordinals. 
These are used in Sections \ref{subsection L-levels}-\ref{subsection role of parameters} 
to prove the characterizations of $\delta$ stated in Theorem \ref{main theorem} 
and the lower bounds for $\delta$ in Theorem \ref{bound on delta}.

For reading this paper, we assume that the reader is familiar with infinite time Turing machines, basic facts about G\"odel's constructible universe and the proof of Welch's submodel characterization of $\lambda$, $\zeta$ and $\Sigma$ from \cite[Theorem 30 \& Corollary 32]{MR2493990}. The latter is used in the proof of Theorem \ref{lambda-zeta-Sigma submodel characterization}.

\subsection{The setting} \label{section setting} 

We briefly introduce the main notions and results related to $\alpha$-ITTMs and refer the reader to \cite{HL, KS, MR3210081, MR2493990} for details. 
We always assume that the tape length $\alpha$ is infinite and multiplicatively closed, i.e. $\mu\cdot \nu<\alpha$ for all $\mu,\nu<\alpha$. It is easy to see that this is equivalent to closure under G\"odel pairing. 
An $\alpha$-ITTM has three tapes of length $\alpha$ for input, working space and output 
and each cell can contain $0$ or $1$. Programs for $T_\alpha$ are just regular Turing machine programs. 
The machine can process a subset of $\alpha$ by representing it on the tapes via its characteristic function. 
Thus we will freely identify a set with its characteristic function. 
The \emph{input tape} carries the subset of $\alpha$ that is given to the machine at the start of the computation, and its content is never changed, while the results of a computation are written on the \emph{output tape}. 
The remaining tape is a \emph{work tape}. 
Furthermore, each tape has a head for reading and writing, all of which move independently of each other. 
It is easy to see that one can equivalently allow any finite number or in fact $\alpha$ many work tapes (using that $\alpha$ is multiplicatively closed). Moreover, the model from \cite{HL} with a single head can simulate our model and is thus equivalent. 

The machine $T_\alpha[e]$ runs along an ordinal time axis.
At successor times, the configuration of the machine is obtained from the preceding one, as usual for a Turing machine, with the extra convention that a head is reset to position $0$ if it is moved to the left from a limit position. At limit times, the content of each cell as well as the head positions are determined as the inferior limits of the sequences of earlier contents of that cell and earlier head positions;
if for some head the inferior limit of the sequence of earlier positions is $\alpha$, then it is reset to $0$. 

A \emph{$\hat{T}_\alpha$-program} computes relative to a finite parameter subset $p$ of $\alpha$. This is given to the program by writing the characteristic function of $p$ to one of the work tapes before the computation starts. 
As we will only be concerned with the case that $\alpha$ is closed under the G\"odel pairing function and the function's restriction to $\alpha$ is easily seen to be computable by an $\alpha$-ITTM, we can assume that parameters are single ordinals below $\alpha$.

We now turn to various notions of writability from \cite{HL}. 
A subset $x$ of $\alpha$ is called \emph{$T_\alpha$-writable} if there is a $T_\alpha$-program $P$ that halts with $x$ on the output tape when the initial input is empty, i.e., all cells contain $0$. 
Moreover, $x$ is called \emph{eventually $T_\alpha$-writable} if there is a $T_\alpha$-program $P$ such that the output tape will have the contents $x$ and never change again from some point on, if the initial input is empty, although the contents of other tapes might change. 
Finally, $x$ is called \emph{accidentally $T_\alpha$-writable} if there is an $T_\alpha$-program such that $x$ appears as the content of the output tape at some time of the computation with empty input. 
Analogous to \cite[Theorem 3.8]{HL}, these three notions of writability are distinct (see Lemma \ref{lambda zeta Sigma with and without parameters}). 

As for \ITTM s, an ordinal is called \emph{$T_\alpha$-clockable} if it is the halting time of a $T_\alpha$-program  with input $\vec{0}$. 

The above notions are defined for $\hat{T}_\alpha$ in an analogous way. 

As for Turing machines, there is a \emph{universal $T_\alpha$-program} $\UU_\alpha$ that simulates all computations with empty input. This can be obtained by dividing the work and output tapes into infinitely many tapes of the same length. 
Note that any $\hat{T}_\alpha$-program can be simulated by a $T_\alpha$-program that considers all possible parameters. 
Thus $U_\alpha$ accidentally writes every $\hat{T}_\alpha$-accidentally writable subset of $\alpha$. 

To compare the writability strength of these machines for different ordinals, we often consider $T_\alpha$-writable subsets $x$ of some ordinal $\xi\leq\alpha$. Naively, we could just write $x$ to the initial segment of length $\xi$ of the output tape and leave the rest empty, but then we could no longer distinguish between $x$ as a subset of $\xi$ and as a subset of $\alpha$. Therefore, we introduce the following notion. 
A subset $x$ of $\xi$ is called \emph{$T_\alpha$-writable as a subset of $\xi$} if there is a $T_\alpha$-program with empty input that halts with the characteristic function of $x$ on the output tape, and if $\xi<\alpha$, then the head is in position $\xi$ at the end of the computation. 
Similarly, we call $x$  \emph{eventually writable as a subset of $\xi$} if the contents of the output tape eventually stabilizes at the characteristic function of $x$ and the head on the output tape eventually stabilizes at $\xi$. 
For any $\xi\leq \min\{\alpha,\beta\}$, we say that $T_\alpha$ has strictly greater writability strength than $T_\beta$ with respect to subsets of $\xi$ if every subset of $\xi$ that is $T_\beta$-writable as a subset of $\xi$ is also $T_\alpha$-writable as a subset of $\xi$, but not conversely. 

Moreover, we frequently use codes for ordinals. 
An \emph{$\alpha$-code} is a subset of $\alpha$ interpreted as a binary relation $\in_\alpha$ on $\alpha$ 
via G\"odel pairing such that $(\alpha,\in_\alpha)$ is well-founded and extensional. 
This structure is isomorphic to a transitive set. The coded set is defined as the image of $0$ in the transitive collapse. 
We work with the image of $0$ instead of the whole set, since the former allows us to code arbitrary sets (this is necessary in Lemmas \ref{equivalent definitions of mugamma} and \ref{equivalent definitions of mu<}), while the latter would only yield codes for transitive sets. 

We further call an ordinal \emph{$T_\alpha$-writable}, \emph{$T_\alpha$-eventually writable} or \emph{$T_\alpha$-accidentally writable} if it has an $\alpha$-code with the corresponding property.\footnote{
Note that the present terminology differs from that of~\cite{ MR3210081}, in which \emph{$T_\alpha$-writability} and \emph{$T_\alpha$-eventual writability} referred to $\omega$-length binary output sequences (as in~\cite{HL}), and $\omega$-codes rather than $\alpha$-codes represented ordinals (and only countable ordinals were considered).  
Results from there need not hold for the current sense of $T_\alpha$-writability, $T_\alpha$-eventual writability, etc. }
The analogous notions for $\hat{T}_\alpha$ are defined similarly. 
Note that one could similarly talk about writability for arbitrary sets, 
and it is easy to see that for sets of ordinals, this would agree with the original definition of writability. 
However, for clarity we will only use this terminology for ordinals. 




\section{Writable and clockable ordinals} \label{section writable and clockable ordinals} 

In this section, we study variants of writability for $T_\alpha$ and $\hat{T}_\alpha$, the associated ordinals, their characterizations and connections with clockable ordinals.  


The ordinals $\lambda$, $\zeta$ and $\Sigma$, which play an important role in the study of infinite time Turing machines, have analogues for $\alpha$-tape machines. We define $\hat{\lambda}_\alpha$, $\hat{\zeta}_\alpha$, $\hat{\Sigma}_\alpha$ and $\lambda_\alpha$, $\zeta_\alpha$, $\Sigma_\alpha$ as the suprema of the $T_\alpha$-writable, $T_\alpha$-eventually writable and $T_\alpha$-accidentally writable ordinals (with respect to $\alpha$-codes) with and without ordinal parameters, respectively. 

We will prove some basic properties of these ordinals. 
Similar to the case of \ITTM s in \cite[Section 2]{MR2493990}, an $\alpha$-word (i.e. an $\alpha$-length bit sequence) is an element of $L_{\hat{\lambda}_\alpha}$, $L_{\hat{\zeta}_\alpha}$ or $L_{\hat{\Sigma}_\alpha}$ if and only if it is $\hat{T}_\alpha$-writable, $\hat{T}_\alpha$-eventually writable or $\hat{T}_\alpha$-accidentally writable, respectively. 
For $\hat{\lambda}_\alpha$, this follows immediately from Lemma \ref{sup of clockables equals sup of writables with parameters}, for $\hat{\zeta}_\alpha$ from Lemma \ref{sup of stabilization times equals sup of writables with parameters} and for $\hat{\Sigma}_\alpha$ from Lemma \ref{loop} below. 

Moreover, an $\alpha$-word is an element of $L_{\lambda_\alpha}$, $L_{\zeta_\alpha}$ or $L_{\Sigma_\alpha}$ if and only if it is contained as an element in some set with a $T_\alpha$-writable, $T_\alpha$-eventually writable or $T_\alpha$-accidentally writable code, respectively. 
For $\lambda_\alpha$ this follows from Lemma \ref{sup of clockables equals sup of writables} and for $\zeta_\alpha$ from Lemma \ref{sup of stabilization times equals sup of writables without parameters}. The claim for $\Sigma_\alpha$ follows from the previous one about $\hat{\Sigma}_\alpha$ by Lemma \ref{lambda zeta Sigma with and without parameters}. 

Given the previous characterization, the reader might wonder whether all elements of $L_{\lambda_\alpha}$ are necessarily $T_\alpha$-writable. 
This holds if and only if $T_\alpha$ reaches all its cells: 
if every ordinal below $\alpha$ is $T_\alpha$-writable, then one can reach any cell via a program that searches for an isomorphism with an $\alpha$-code for the given ordinal, and the converse is easy to see. 

We will frequently use the fact that for any multiplicatively closed $\xi$, any $\gamma$ with a $T_\alpha$-writable $\xi$-code, $L_\gamma$ also has a $T_\alpha$-writable $\xi$-code, and the same holds for $\hat{T}_\alpha$-writable codes. 
To see this, one partitions $\xi$ into $\gamma$ many pieces with order type $\xi$ and successively writes $\xi$-codes for $L_\mu$ onto the $\mu$th piece for all $\mu<\gamma$. 

The next lemma is used to prove that $\hat{\lambda}_\alpha$ equals the supremum of $\hat{T}_\alpha$-clockable ordinals. 
It shows that any $\hat{T}_\alpha$-program that does not halt on input $\vec 0$ runs into a loop between $\hat{\zeta}_\alpha$ and $\hat{\Sigma}_\alpha$, as for standard \ITTM s. 

\begin{lem} \label{loop} 
On input $\vec 0$, any $\hat{T}_\alpha$-program either 
halts before time $\hat{\zeta}_\alpha$ or runs into an ever-repeating loop in which the configuration at time $\hat{\zeta}_\alpha$ is the same as that of time $\hat{\Sigma}_\alpha$.
\end{lem}
\begin{proof} 
We refer the reader to the proof of this fact for \ITTM s \cite[Section 2]{MR1734198} and only sketch the changes that are necessary to adapt it to $\alpha$-\ITTM s. 
Since ordinal parameters are allowed in the definitions of $\hat{\lambda}_\alpha$, $\hat{\zeta}_\alpha$ and $\hat{\Sigma}_\alpha$, it is sufficient to prove that the limit behaviour in each cell is the same when the time approaches $\hat{\zeta}_\alpha$ and $\hat{\Sigma}_\alpha$. In other words, if the contents of the $\xi$th cell converges when the time appoaches $\hat{\zeta}_\alpha$, then it converges to the same value at $\hat{\Sigma}_\alpha$, otherwise it diverges at $\hat{\Sigma}_\alpha$. 
We need ordinal parameters, since an \ITTM\ with parameter $\xi$ is used to observe the $\xi$th cell. 

The difference to the setting of \ITTM s is that here the head doesn't move to the first cell at every limit time. We want to show that for any computation of $T_\alpha$, the head position at time $\hat{\zeta}_\alpha$ is equal to the head position at time $\hat{\Sigma}_\alpha$. 
To adapt the proof, we define a program that simulates the given machine, and writes the current head position on an additional tape by writing $1$ in every cell that precedes the head position and $0$ everywhere else.
At every limit time, the inferior limit of the head positions is calculated and the contents of the remaining cells are deleted. 
Now the proof for \ITTM s shows that the tape contents for the simulation are identical at the times $\hat{\zeta}_\alpha$ and $\hat{\Sigma}_\alpha$ and thus the head positions are also equal for the original program. 
\end{proof} 

Note that the version of the previous lemma for $\zeta_\alpha$ and $\Sigma_\alpha$ fails if $\zeta_\alpha<\hat{\zeta}_\alpha$: 
a universal $T_\alpha$-program (see Section \ref{section setting}) simulates all $\hat{T}_\alpha$-programs and thus its first ever-repeating loop begins at $\hat{\zeta}_\alpha$. 
Moreover, $\zeta_\alpha<\hat{\zeta}_\alpha$ is possible by Section \ref{subsection role of parameters} below. 

The fact that the suprema of writable and clockable ordinals are equal \cite[Theorem 1.1]{MR1734198} easily generalizes as follows to the setting with ordinal parameters. 

\begin{lem} \label{sup of clockables equals sup of writables with parameters} 
$\hat{\lambda}_\alpha$ equals the strict supremum of $\hat{T}_\alpha$-clockable ordinals. 
\end{lem} 
\begin{proof} 
Let $\xi$ denote the supremum of $\hat{T}_\alpha$-clockable ordinals. 

To show $\hat{\lambda}_\alpha\leq \xi$, take any $T_\alpha$-writable ordinal $\beta$. 
The following program halts after at least $\beta$ steps. 
The program writes an $\alpha$-code for $\beta$, counts through the code by successively deleting the next remaining element and halts when all elements are deleted. 

To show $\hat{\lambda}_\alpha\geq \xi$, take any $\hat{T}_\alpha$-clockable ordinal $\beta$. 
Let $P$ be a $\hat{T}_\alpha$-program that halts at time $\beta$. 
By Lemma \ref{loop}, $\beta<\hat{\zeta}_\alpha$. 
Thus there is an eventually $\hat{T}_\alpha$-writable ordinal $\gamma>\beta$. 

Consider the following $\hat{T}_\alpha$-program. 
The program writes each version $\mu$ of $\gamma$ and runs $P$ up to time $\mu$. Whenever $\mu$ changes, begin a new simulation. 
It is clear that this will halt when $\mu\geq\beta$. When this happens, output an $\alpha$-code for $\mu$. Thus $\mu$ is $\hat{T}_\alpha$-writable. 

Since $\hat{\lambda}_\alpha$ is itself not $\hat{T}_\alpha$-writable, the previous argument shows that the supremum is strict. 
\end{proof} 

\begin{lem} \label{sup of stabilization times equals sup of writables with parameters} 
$\hat{\zeta}_\alpha$ equals the strict supremum of $\hat{T}_\alpha$-stabilization times of the tape contents. 
\end{lem} 
\begin{proof} 
Let $\gamma$ denote the supremum of $\hat{T}_\alpha$-stabilization times of the tape contents. 

To show that $\hat{\zeta}_\alpha\leq\gamma$, suppose that $P$ eventually writes $\xi$. We consider a program $Q$ that simulates $P$ and additionally sets a flag. It is set to $0$ when the output of $P$ changes and to $1$ once we have counted through $\mu$, if the current output of $P$ codes an ordinal $\mu$. Then $Q$'s stabilization time is at least $\xi$. 

By Lemma \ref{loop}, $\gamma\leq \hat{\zeta}_\alpha$. 
To show that the supremum is strict, assume that some program $P$ stabilizes exactly at time $\hat{\zeta}_\alpha$. 
We search via the universal $T_\alpha$-program for (an accidentally writable code for) some $\mu$ such that the output of $P$ eventually stabilizes before $\mu$. 
Then $\mu$ is eventually $\hat{T}_\alpha$-writable, contradicting the fact that $\mu>\hat{\zeta}_\alpha$. 
\end{proof} 

The next result describes the basic relations between the ordinals associated to $T_\alpha$ and $\hat{T}_\alpha$. 
 
\begin{lem} \label{lambda zeta Sigma with and without parameters} \hfill 
\begin{enumerate} 
\item \label{lambda zeta Sigma with and without parameters 1}
$\hat{\lambda}_\alpha$ is $T_\alpha$-eventually writable. 
\item \label{lambda zeta Sigma with and without parameters 2}
$\hat{\zeta}_\alpha$ is $T_\alpha$-accidentally writable. 
\item \label{lambda zeta Sigma with and without parameters 3}
$\hat{\Sigma}_\alpha=\Sigma_\alpha$. 
\end{enumerate} 
Therefore 
$\lambda_\alpha\leq \hat{\lambda}_\alpha<\zeta_\alpha\leq \hat{\zeta}_\alpha<\Sigma_\alpha=\hat{\Sigma}_\alpha$. 
\end{lem} 
\begin{proof} 
To show \ref{lambda zeta Sigma with and without parameters 1}, we simulate all $\hat{T}_\alpha$-programs, beginning with the first step of each computation and proceeding with one step of each program at a time. 
This is done by partitioning the tape into $\alpha$ many tapes of length $\alpha$. 
For each $i<\alpha$, we define $\gamma_{i,j}$ as follows. If the $i\mathrm{th}$ $\hat{T}_\alpha$-program $P_i$ halts in step $j$ of the run of $P_i$ with output a (code for an) ordinal $\gamma$, let $\gamma_{i,j}=\gamma$. Otherwise let $\gamma_{i,j}=0$. 
Let further $\gamma_j=\sum_{i<\alpha} \gamma_{i,j}$. 
The output of our algorithm is set to the value $\gamma[j]=\sum_{i\leq j} \gamma_i$ once the $j\mathrm{th}$ step of each program is completed. 

In step $\hat{\lambda}_\alpha$ of the simulation, all steps $j<\hat{\lambda}_\alpha$ of each $\hat{T}_\alpha$-program are completed. 
Moreover, each $\hat{T}_\alpha$-program has either already halted or diverges by Lemma \ref{sup of clockables equals sup of writables with parameters}. 
Hence the output of the simulation takes the constant value $\gamma=\sum_{j<\hat{\lambda}_\alpha} \gamma_j=\sup_{j<\hat{\lambda}_\alpha} \gamma[j]$ from step $\hat{\lambda}_\alpha$ onwards. 

It remains to show that $\gamma=\hat{\lambda}_\alpha$. 
To see that $\gamma\leq\hat{\lambda}_\alpha$, note that $\gamma[j]<\hat{\lambda}_\alpha$ for all $j<\hat{\lambda}_\alpha$, since $\gamma[j]$ is $\hat{T}_\alpha$-writable. 
To see that $\gamma\geq\hat{\lambda}_\alpha$, note that every $\hat{T}_\alpha$-writable ordinal is of the form $\gamma_{i,j}$ for some $i<\alpha$ and $j<\hat{\lambda}_\alpha$ and $\gamma_{i,j}\leq\gamma_j\leq\gamma$. 

The proof of \ref{lambda zeta Sigma with and without parameters 2} is similar. 
We simulate all $\hat{T}_\alpha$-programs as above. For each $i<\alpha$, let $\gamma_{i,j,k}$ denote the output of the $i{\mathrm{th}}$ program $P_i$ in step $k$ of the run of $P_i$, if this codes an ordinal, is constant in the interval $[j,k)$ and $j$ is minimal with this property. 
Let $\gamma_{i,j,k}=0$ otherwise. 
Let further $\gamma_{j,k}=\sum_{i<\alpha}\gamma_{i,j,k}$. 
The algorithm's output is set to $\gamma[k]=\sum_{j<k} \gamma_{j,k}$ after the $j\mathrm{th}$ steps of each program are completed for all $j<k$. 
Now let $k=\hat{\zeta}_\alpha$. 

To see that $\gamma[k]\leq\hat{\zeta}_\alpha$, note that $\sum_{j<l}\gamma_{j,k}$ is $\hat{T}_\alpha$-eventually writable for all $l<k$ by Lemma \ref{loop}. 
To see that $\gamma[k]\geq\hat{\zeta}_\alpha$, note that every $\hat{T}_\alpha$-eventually writable ordinal is of the form $\gamma_{i,j,k}$ for some $i<\alpha$ and $j<k$ by Lemma \ref{sup of stabilization times equals sup of writables with parameters}.\footnote{We did not specify at which time of the simulation the output equals $\gamma[k]$. It can be arranged that this happens at time $k$ for $k=\hat{\zeta}_\alpha$. }   

Finally, \ref{lambda zeta Sigma with and without parameters 3} follows from the fact that any $\hat{T}_\alpha$-accidentally writable subset of $\alpha$ is $T_\alpha$-accidentally writable. This was already shown in Section \ref{section setting}. 
\end{proof} 

\noindent
We will see that $\lambda_\alpha<\hat{\lambda}_\alpha$ and $\zeta_\alpha<\hat{\zeta}_\alpha$ for some $\alpha$ in Section \ref{subsection role of parameters}. 

We can now prove a version of Theorem \ref{sup of clockables equals sup of writables with parameters} without parameters. 

\begin{lem} \label{sup of clockables equals sup of writables} 
$\lambda_\alpha$ equals the supremum of $T_\alpha$-clockable ordinals. 
\end{lem} 
\begin{proof} 
Let $\gamma$ denote the supremum of $T_\alpha$-clockable ordinals. 

It is easy to see that $\lambda_\alpha\leq\gamma$ (as for $\hat{\lambda}_\alpha$). 

To see that $\gamma\leq\lambda_\alpha$, take any $T_\alpha$-clockable ordinal $\xi$. Then $\xi<\hat{\lambda}_\alpha$ by Lemma \ref{sup of clockables equals sup of writables with parameters} and thus $\xi<\zeta_\alpha$ by Lemma \ref{lambda zeta Sigma with and without parameters}. 
Fix a $T_\alpha$-program $P$ halting at time $\xi$ and a $T_\alpha$-program $Q$ that eventually writes some $\mu\geq\xi$. For each ordinal $\nu$ output by $Q$, we simulate $P$ up to time $\nu$ and output $\nu$ if $P$ halts. The output is a $T_\alpha$-writable ordinal $\nu\geq\xi$. Thus $\xi<\lambda_\alpha$. 
\end{proof} 

Similarly as in  Lemma \ref{sup of stabilization times equals sup of writables with parameters} for $\hat{\zeta}_\alpha$, we obtain the following version without parameters.  

\begin{lem} \label{sup of stabilization times equals sup of writables without parameters} 
$\zeta_\alpha$ equals the strict supremum of $T_\alpha$-stabilization times of the tape contents. 
\end{lem} 
\begin{proof} 
Let $\gamma$ denote the supremum of $T_\alpha$-stabilization times of the tape contents. 

It is easy to see that $\zeta_\alpha\leq\gamma$ (as for $\hat{\zeta}_\alpha$). 

We have $\gamma\leq \hat{\zeta}_\alpha$ by Lemma \ref{loop} and hence $\gamma<\Sigma_\alpha$ by Lemma \ref{lambda zeta Sigma with and without parameters}. 
The next argument for the inequality $\gamma\leq \zeta_\alpha$ and for the fact that this is a strict supremum is virtually the same as in Lemma \ref{sup of stabilization times equals sup of writables with parameters}. 
Suppose that some program $P$ stabilizes at a time $\eta\geq\zeta_\alpha$. 
Since $\eta\leq\gamma<\Sigma_\alpha$, $\eta$ is accidentally $T_\alpha$-writable. 
We search via the universal $T_\alpha$-program for (an accidentally writable code for) the least $\mu$ such that the output of $P$ eventually stabilizes at time $\mu$. 
Then $\mu=\eta$ is eventually $T_\alpha$-writable, contradicting the fact that $\eta\geq\zeta_\alpha$. 
\end{proof}

We next show that $\hat{\lambda}_\alpha$ is admissible and $\hat{\zeta}_\alpha$ is $\Sigma_2$-regular. 
We first fix some notation. 
Given a class $\Sigma$ of formulas, an ordinal $\gamma$ is called \emph{$\Sigma$-regular} if for no $\beta<\gamma$, there is a cofinal function $f\colon \beta\rightarrow \gamma$ that is $\Sigma$-definable over $L_\gamma$ from parameters in $L_\gamma$. 
Moreover, $\Sigma_1$-regular ordinals are called \emph{admissible}. 
To show that $\hat{\lambda}_\alpha$ is admissible, we need the following lemma (which must be folklore). 
To state the lemma, recall that \emph{$\Sigma$-collection} states that for any $\Sigma$-formula $\varphi(x,y)$ and set $A$ with $\forall x\in A\ \exists y\ \varphi(x,y)$, there is a set $B$ with $\forall x\in A\ \exists y\in B \ \varphi(x,y)$. 

\begin{lem} \label{equivalence of admissible and collection} 
Let $\gamma\in\Ord$ and $n\in\omega$. The following statements are equivalent: 
\begin{enumerate}[(a)] 
\item \label{equivalence of admissible and collection 1} 
$\gamma$ is $\Sigma_{n+1}$-regular.\footnote{A formula is called $\Sigma_0$ if it contains only bounded quantifiers, 
$\Sigma_{n+1}$ if it logically equivalent to a formula of the form $\exists x_0,\dots,x_n \varphi$, where $\varphi$ is $\Pi_n$, and $\Pi_n$ if it is logically equivalent to a formula of the form $\neg \varphi$, where $\varphi$ is $\Sigma_n$. 
It follows that these classes of formulas are closed under the connectives $\wedge$ and $\vee$. }
\item \label{equivalence of admissible and collection 2} 
$L_\gamma\models \Pi_n$-{\rm collection}. 
\item \label{equivalence of admissible and collection 3} 
$L_\gamma\models \Sigma_{n+1}$-{\rm collection}. 
\end{enumerate} 
\end{lem} 
\begin{proof} 
Assume that \ref{equivalence of admissible and collection 1} holds. 
To show \ref{equivalence of admissible and collection 2}, 
take a $\Pi_n$-formula $\varphi(x,y,z)$ and $A,B \in L_\gamma$ with $L_\gamma\models \forall x\in A\ \exists y\ \varphi(x,y,B)$. 
Let further $f_A \colon \gamma_A\rightarrow A$ denote the order-preserving enumeration of $A$ with respect to $\leq_L$ 
and $f\colon \gamma_A \rightarrow \gamma$ the function with $f(\alpha)$ equal to the least $\beta<\gamma$ with $L_\gamma\models \exists y\in L_\beta\ \varphi(f_A(\alpha),y,B)$. 
Then $f_A$ is $\Delta_1$-definable over $L_\gamma$  from $A$. 

We claim that $f$ is $\Delta_{n+1}$-definable over $L_\gamma$ from $A,B$. 
Note that it follows from $\Sigma_k$-regularity by induction on $i\leq k$ that $\Sigma_i$- and $\Pi_i$-formulas are closed under bounded quantification. 
Thus $\exists y\in L_\beta\ \varphi(f_A(\alpha),y,B)$ is (in $L_\gamma$) equivalent to a $\Pi_n$-formula and the function sending $\alpha$ to $L_{f(\alpha)}$ is definable by the conjunction of a $\Sigma_n$ and a $\Pi_n$-formula, so it is $\Delta^1_{n+1}$-definable over $L_\gamma$ from $A,B$. 
It follows that $f$ is $\Delta^1_{n+1}$-definable over $L_\gamma$ from $A,B$. 

By \ref{equivalence of admissible and collection 1}, $\ran{f}$ is bounded by some $\beta<\gamma$. 
Thus $L_\beta$ witnesses $\Pi_0$-collection for $\varphi(x,y,z)$ and $A,B$. 
Therefore \ref{equivalence of admissible and collection 2} holds. 

It is easy to see that \ref{equivalence of admissible and collection 2} implies \ref{equivalence of admissible and collection 3} and \ref{equivalence of admissible and collection 3} implies \ref{equivalence of admissible and collection 1}. 
\end{proof} 
\noindent
The next result is analogous to the fact that $\lambda$ is admissible \cite[Corollary 8.2]{HL}. 



\begin{lem}\footnote{An anonymous referee asked whether this also holds for $\lambda_\alpha$ and $\zeta_\alpha$. }\label{lambdaalpha is admissible}
\begin{enumerate} 
\item \label{lambda_alpha is admissible 1} 
$\hat{\lambda}_\alpha$ is admissible. 
\item \label{lambda_alpha is admissible 2} 
$\hat{\zeta}_\alpha$ is $\Sigma_2$-regular. 
\end{enumerate} 
\end{lem} 
\begin{proof} 
To show that $\hat{\lambda}_\alpha$ is admissible, it is sufficient to show that $L_{\hat{\lambda}_\alpha}$ is a model of $\Pi_0$-collection by Lemma \ref{equivalence of admissible and collection}. 
To this end, take any $\Pi_0$-formula $\varphi(x,y,z)$ and $A,B\in L_{\hat{\lambda}_\alpha}$ with $L_{\hat{\lambda}_\alpha}\models \forall x\in A\ \exists y\ \varphi(x,y,B)$. 
Thus $A$ has a $\hat{T}_\alpha$-writable code. 
We generate outputs $\gamma$ via the universal $T_\alpha$-program. 
When $A\in L_\gamma$ and $L_\gamma\models \forall x\in A\ \exists y\ \varphi(x,y,B)$, output $\gamma$ and halt. 
This program will halt since $\hat{\lambda}_\alpha$ is $T_\alpha$-accidentally writable by Lemma \ref{lambda zeta Sigma with and without parameters}, thus producing some $\gamma<\hat{\lambda}_\alpha$. Hence $\Pi_0$-collection holds in $L_{\hat{\lambda}_\alpha}$. 

The proof of \ref{lambda_alpha is admissible 2} is similar. 
By Lemma \ref{equivalence of admissible and collection}, it is sufficient to show $\Pi_1$-collection in $L_{\hat{\zeta}_\alpha}$. 
To see this, take a $\Pi_1$-formula $\varphi(x,y,z)$ and $A,B\in L_{\hat{\zeta}_\alpha}$ with $L_{\hat{\zeta}_\alpha}\models \forall x\in A\ \exists y\ \varphi(x,y,B)$. 
Thus $A, B\in L_{\hat{\zeta}_\alpha}$ have $\hat{T}_\alpha$-eventually writable codes. 
Consider the following $\hat{T}_\alpha$-program. 
For the current versions of $A$ and $B$, we search for (a code for) an ordinal $\gamma$ via the universal $T_\alpha$-program $U_\alpha$ and simultaneously for each $x\in A$ for some $y\in L_\gamma$ with $\varphi(x,y,B)$. 
More precisely, we implement the following (simultaneous) subroutines for all $x\in A$. 
Take $x\in A$ and $y\in L_\gamma$ as the current candidate for $\varphi(x,y,B)$. 
We run a search for counterexamples to $\varphi(x,y,B)$ via $U_\alpha$; if no counterexample is found, then we keep $y$, but discard it otherwise. 
If the subroutines eventually stabilize for all $x$, then $\gamma$ is the eventual output. 
On the other hand, there might be some $x\in A$ such that all its candidates are discarded at some time; we then move on to $\gamma+1$. 
Clearly there is an eventual output $\gamma<\hat{\zeta}_\alpha$. 
Hence $\Pi_1$-collection holds in $L_{\hat{\zeta}_\alpha}$. 
\end{proof} 
\noindent
We will further use the next variant of the submodel characterisation of $\lambda$, $\zeta$ and $\Sigma$. 
We say that a tuple $(\alpha_0,\dots,\alpha_n)$ is \emph{least} with a certain property if $\forall i\leq n\ \alpha_i\leq\beta_i$ for any other such tuple $(\beta_0,\dots,\beta_n)$. 

\begin{thm} \label{lambda-zeta-Sigma submodel characterization}
$(\hat{\lambda}_{\alpha},\hat{\zeta}_{\alpha},\hat{\Sigma}_{\alpha})$ is 
the least triple $(\mu,\nu,\xi)$ with $\alpha<\mu<\nu<\xi$ and $L_\mu \prec_{\Sigma_{1}} L_\nu\prec_{\Sigma_{2 }}L_\xi$. 
\end{thm}
\begin{proof} 
The proof of $L_{\hat{\lambda}_{\alpha}}\prec_{\Sigma_{1}} L_{\hat{\zeta}_{\alpha}}\prec_{\Sigma_{2}}L_{\hat{\Sigma}_{\alpha}}$ is virtually the same as for $(\lambda,\zeta,\Sigma)$ in \cite[Corollary 32]{MR2493990}. 

The proof of minimality of $\hat{\zeta}_\alpha$ in \cite[Theorem 30]{MR2493990} for $\alpha=\omega$ adapts to this setting. 
We briefly discuss the crucial role of parameters in our version. 
The distinction between computations with and without parameters is not visible in Welch's proof, as finite parameters are always writable. 
First, to show that the content of a tape cell stabilizes at time $\hat{\zeta}_{\alpha}$ if and only if it stabilizes at time $\hat{\Sigma}_{\alpha}$, it is necessary to let the machine check  the evolution of the contents of each cell separately for each cell as in Lemma \ref{loop}. 
This is clearly possible for the $\xi$th cell if $\xi$ is given as a parameter. 
Second, it is frequently needed that any element $x$ of a set $y$ with a $\hat{T}_\alpha$-writable code has itself a $\hat{T}_\alpha$-writable code. 
This need not be true for $T_\alpha$, as $x$ might correspond to an ordinal in the code for $y$ that is not $T_{\alpha}$-reachable. 
However, the statement for $\hat{T}_\alpha$-writable codes and its analogue for $\hat{T}_\alpha$-eventually writable codes clearly hold for $\hat{T}_\alpha$. 
Finally, for our machines the read-write-head is no longer reset to $0$ at all limit times. 
In the $\omega$-case, this is used to show that the snapshots at times $\zeta$ and $\Sigma$ agree.
But this issue has already been dealt with in the proof of Lemma \ref{loop}. 

To see that $\hat{\lambda}_\alpha$ is also minimal, take a triple $(\mu,\nu,\xi)$ as above. Since $\nu\geq\hat{\zeta}_\alpha$, every halting $\hat{T}_\alpha$-program halts before $\mu$ and hence $\mu\geq\hat{\lambda}_\alpha$. 
To finally see that $\hat{\Sigma}_\alpha$ is minimal, suppose that $(\mu,\nu,\xi)$ is a triple with $\xi<\hat{\Sigma}_\alpha$. Since $L_{\hat{\lambda}_\alpha}\prec_{\Sigma_1}L_{\hat{\Sigma}_\alpha}$ there is such a triple below $\hat{\lambda}_\alpha$, but this contradicts the fact that $\nu\geq\hat{\zeta}_\alpha$. 
\end{proof} 

Is there a version of the previous result for $(\lambda_\alpha,\zeta_\alpha,\Sigma_\alpha)$? This was asked by one of the referees of this paper. 
For this triple, it is natural to consider the class $\Sigma_n^{(\alpha)}$ of formulas with parameter $\alpha$, and in fact $L_{\lambda_\alpha}\prec_{\Sigma_1^{(\alpha)}}L_{\zeta_\alpha}\prec_{\Sigma_2^{(\alpha)}}L_{\Sigma_\alpha}$ remains valid. 
To see that $(\lambda_\alpha,\zeta_\alpha,\Sigma_\alpha)$ is not necessarily the least such triple, suppose that $\zeta_\alpha<\hat{\zeta}_\alpha$ (this is possible by Lemma \ref{role of parameters} below). 
We have $L_{\zeta_\alpha}\prec_{\Sigma_2^{(\alpha)}}L_{\hat{\zeta}_\alpha}$, since $L_{\zeta_\alpha}\prec_{\Sigma_2^{(\alpha)}}L_{\Sigma_\alpha}$ and $L_{\hat{\zeta}_\alpha}\prec_{\Sigma_2^{(\alpha)}}L_{\hat{\Sigma}_\alpha}$. 
Hence $(\lambda_{\alpha},\zeta_{\alpha},\hat{\zeta}_{\alpha})$ is a triple with the required property, but $\hat{\zeta}_\alpha<\Sigma_\alpha$ by Lemma \ref{lambda zeta Sigma with and without parameters}.

\section{
Writability strength, reachability and $L$-levels} \label{section writability strength, reachability and L-levels} 

\subsection{Local cardinals} \label{subsection L-levels}

We characterize $\delta$ by connecting properties of levels of the constructible universe with writability strength. 
To aid this, we begin with some elementary observations. 

Let $\Card_*$ denote the set of ordinals $\alpha>\omega$ that are cardinals in $L_{\lambda_\alpha}$. The next observation states some properties of this set. 
We will see in Section \ref{subsection reachable cells} that $\delta=\min(\Card_*)$. 

\begin{obs} \label{Obs 1} 
Suppose that $\kappa$ is an uncountable cardinal. 
\begin{enumerate} 
\item \label{Obs 1.1}
$\Card_*$ is unbounded in $\kappa$.  
\item \label{Obs 1.2} 
For any $\alpha<\kappa$, there is a sequence of length $\alpha$ of successive multiplicatively closed ordinals below $\kappa$ that is disjoint from $\Card_*$.\footnote{We defined $\lambda_\alpha$ only for multiplicatively closed ordinals, so the remaining ordinals are by definition not elements of $\Card_*$. } 
\end{enumerate} 
\end{obs} 

\begin{proof} 
It is sufficient to prove \ref{Obs 1.1} assuming that $\kappa$ is regular. Take any $\xi<\kappa$ and let $\pi\colon h^{L_{(\kappa^+)^L}}(\xi+1)\rightarrow L_\beta$ denote the transitive collapse. Then $\alpha=\pi(\kappa)$ is a cardinal in $L_\beta$ and $\xi<\alpha<\lambda_\alpha<\beta$. 
Hence $\alpha\in\Card_*$. 

For \ref{Obs 1.2}, we again take a regular $\kappa$. 
Let $\beta=\omega^{\omega^\alpha}$, $1\leq\eta<\alpha$ and $\gamma=\omega^{\omega^{\alpha+\eta}}$ the $\eta$th multiplicatively closed ordinal above $\beta$. 
Note that the inductive definition of multiplication can be carried out in $L_\gamma$, since $\gamma$ is multiplicatively closed. 
It follows that there is a definable (over $L_\gamma$) surjection from $\eta$ onto the set of multiplicatively closed ordinals between $\beta$ and $\gamma$. 

Moreover, it is easy to see that there are uniformly in $n\in\omega$ definable (over $L_\gamma$) functions sending ordinals $\theta$ to surjections $\theta\rightarrow\theta^n$. From these, we obtain a function sending $\theta$ to a surjection $\theta\rightarrow\theta^\omega$, the least multiplicatively closed ordinal above $\theta$. 

Using the previous functions, one easily obtains a definable (over $L_\gamma$) surjection from $\beta$ onto $\gamma$. 
In particular, $\gamma$ is not a cardinal in $L_{\lambda_\gamma}$. 
\end{proof} 

\noindent
By Observation~\ref{Obs 1}, the following ordinals are well-defined. 

\begin{defi} \label{definition of mu} \hfill 
\begin{enumerate}[(a)] 
\item 
For any ordinal $\xi$, let $\mu_\xi$ be the least $\alpha$ with $L_{\lambda_\alpha}\models |\alpha|>\xi$. 
\item 
Let $\mu_*=\min(\Card_*)$. 
\end{enumerate}
\end{defi}

\noindent
The next lemma shows that the ordinals in the previous definition are equal for $\xi=\omega$. 

\begin{lem} \label{thetaomega equals theta<}
$\mu_*=\mu_\omega$. 
\end{lem} 
\begin{proof} 
$\mu_\omega\leq \mu_*$ is clear. Assume towards a contradiction that $\mu_\omega<\mu_*$. 
By the definition of $\mu_*$, there is a surjection $f\colon \xi \rightarrow \mu_\omega$ in $L_{\lambda_{\mu_\omega}}$ for some $\xi<\mu_\omega$. 
Consider a $T_{\mu_\omega}$-program that searches for such a $\xi<\mu_\omega$ and a surjection $f\colon \xi \rightarrow \mu_\omega$. 
We fix $\xi$ and $f$ that are found by the program. 

Note that $\xi$ is $T_{\mu_\omega}$-writable, so $T_{\mu_\omega}$ can simulate $T_\xi$. 
Once our search suceeds, we search for a surjection $g\colon \omega\rightarrow \xi$ via a $T_\xi$-program. 
This will also suceed, since such a surjection exists in $L_{\lambda_\xi}$ by the definition of $\mu_\omega$. 
We have produced a surjection $f\circ g\colon \omega\rightarrow \mu_\omega$ in $L_{\mu_\omega}$. 
But this contradicts Definition \ref{definition of mu}. 
\end{proof} 

\begin{lem} \label{equivalent definitions of mugamma} 
The following statements are equivalent: 
\begin{enumerate}[(a)] 
\item \label{condition for mu 1} 
$L_{\lambda_\alpha}\models |\alpha|>\xi$. 
\item \label{condition for mu 2} 
There is no $\alpha$-code for a surjection $f\colon \xi\rightarrow \alpha$ that is $T_\alpha$-writable with $\xi$ as a parameter. 
\item \label{condition for mu 3} 
$L_{\hat{\lambda}_\alpha}\models |\alpha|>\xi$. 
\item \label{condition for mu 4} 
There is no $\hat{T}_\alpha$-writable $\alpha$-code for a surjection $f\colon \xi\rightarrow \alpha$. 
\end{enumerate} 
\end{lem} 
\begin{proof} 
It is easy to see that \ref{condition for mu 1} is equivalent to \ref{condition for mu 2}, \ref{condition for mu 3} to \ref{condition for mu 4} and \ref{condition for mu 4} implies \ref{condition for mu 2}. 
To see that \ref{condition for mu 2} implies \ref{condition for mu 4}, it suffices to write such a code only from $\xi$. 
This can be done by simulating the program for \ref{condition for mu 4} simultaneously for all ordinal parameters and halting when the required code appears. 
\end{proof} 

\noindent
Using the previous lemma, one can observe that $\mu_\xi$ equals the least $\alpha>\xi$ in $\Card_*$ and therefore, the function $\xi\mapsto \mu_\xi$ enumerates the successor elements of $\Card_*$ (i.e. those which are not limits of $\Card_*$).\footnote{We would like to thank an anonymous referee for this observation. } 
To see this, it suffices to show that for $\alpha=\mu_\xi$, we have $L_{\hat{\lambda}_\alpha}\models \xi^+=\alpha$. 
Note that $\nu:=(\xi^+)^{L_{\hat{\lambda}_\alpha}}\leq\alpha$ by the definition of $\mu_\xi$ and Lemma \ref{equivalent definitions of mugamma}. 
Since one can simulate shorter tapes by using ordinal parameters, the function $\gamma\mapsto \hat{\lambda}_\gamma$ is monotone. 
We must then have $\nu=\alpha$, since $\nu<\alpha$ would contradict the minimality of $\alpha$.

We further obtain the next equivalences by virtually the same proof as for Lemma \ref{equivalent definitions of mugamma}.

\begin{lem} \label{equivalent definitions of mu<} 
The following statements are equivalent: 
\begin{enumerate}[(a)] 
\item \label{condition for mu 1} 
$L_{\lambda_\alpha}\models \alpha$ is a cardinal. 
\item \label{condition for mu 2} 
There is no $T_\alpha$-writable $\alpha$-code for a surjection $f\colon \xi\rightarrow \alpha$ for some $\xi<\alpha$. 
\item \label{condition for mu 3} 
$L_{\hat{\lambda}_\alpha}\models \alpha$ is a cardinal. 
\item \label{condition for mu 4} 
There is no $\hat{T}_\alpha$-writable $\alpha$-code for a surjection $f\colon \xi\rightarrow \alpha$ for some $\xi<\alpha$. 
\end{enumerate} 
\end{lem} 

\noindent
We can replace  the surjections in Definition \ref{definition of mu} by cofinal functions. 
This yields results analogous to Lemmas \ref{thetaomega equals theta<}, \ref{equivalent definitions of mugamma} and \ref{equivalent definitions of mu<} with virtually the same proofs, which we do not state explicitly. 

\begin{defi} \label{definition of nu}\hfill 
\begin{enumerate} 
\item 
For any ordinal $\xi$, let $\nu_\xi$ be the least $\alpha$ with $L_{\lambda_\alpha}\models \cof(\alpha)>\xi$. 
\item 
Let $\nu_*$ be the least $\alpha>\omega$ that is regular in $L_{\lambda_\alpha}$. 
\end{enumerate} 
\end{defi} 


\noindent
The previous results yield the next equality. 

\begin{lem}\label{all values are equal} 
$\mu_\omega=\nu_\omega$. 
\end{lem} 
\begin{proof} 
It is clear that $\mu_\omega\leq\nu_\omega$. Assume towards a contradiction that $\mu_\omega<\nu_\omega$. 
We first search for a $T_{\mu_\omega}$-writable cofinal function $f\colon \omega\rightarrow \mu_\omega$ using the analogue to Lemma \ref{equivalent definitions of mugamma} for $\nu_\omega$. 
We then search for a sequence of surjections $f_n\colon \omega\rightarrow f(n)$. 
This will succeed since $f(n)<\mu_\omega$ for all $n\in\omega$ and by the definition of $\mu_\omega$. 
The algorithm yields a $T_{\mu_\omega}$-writable surjection from $\omega$ onto $ \mu_\omega$, contradicting the definition of $\mu_\omega$. 
\end{proof} 

It follows from the combined results in this section that $\mu_*=\nu_*$ is the least $\alpha$ with either of the properties (a) $\alpha$ is uncountable in $L_{\lambda_\alpha}$ (b) $\alpha$ is regular in $L_{\lambda_\alpha}$ or (c) $\alpha$ is a cardinal in $L_{\lambda_\alpha}$. 
This proves that the least ordinals satisfying \ref{main theorem 6}-\ref{main theorem 8} of Theorem \ref{main theorem} are equal. 

To see that these equal the least ordinal with \ref{main theorem 9} of Theorem \ref{main theorem}, first note that for $\hat{\lambda}_\alpha$, this follows from the previous results. 
Moreover, the claim for $\Sigma_\alpha=\hat{\Sigma}_\alpha$ and $\hat{\zeta}_\alpha$ holds, since we have $L_{\hat{\lambda}_\alpha}\prec_{\Sigma_1}L_{\Sigma_\alpha}$ and $L_{\hat{\lambda}_\alpha}\prec_{\Sigma_1} L_{\hat{\zeta}_\alpha}$ by Theorem \ref{lambda-zeta-Sigma submodel characterization}. 
The argument for $\zeta_\alpha$ is analogous to the proofs of Lemmas \ref{equivalent definitions of mugamma} and \ref{equivalent definitions of mu<}. 

\subsection{Reachable cells} \label{subsection reachable cells}

We now give characterizations of $\delta$ via some results in the previous section. 

\begin{prop} \label{delta equals mu} 
$\delta=\mu_*$. 
\end{prop} 
\begin{proof} 
To see that $\delta\leq\mu_*$, it suffices to show that $T_{\mu_*}$ doesn't reach all its cells. We thus assume otherwise. Then there is a well-defined map $f\colon \mu_*\rightarrow \Ord$ that sends each $\alpha<\mu_*$ to the least halting time of a program that halts with its head in the $\alpha$th cell. Since the values are bounded by $\lambda_{\mu_*}\leq\hat{\lambda}_{\mu_*}$ by Lemma \ref{sup of clockables equals sup of writables}, $f$ is $\Sigma_1$-definable over $L_{\hat{\lambda}_{\mu_*}}$. Since $\hat{\lambda}_{\mu_*}$ is admissible by Lemma \ref{lambdaalpha is admissible}, $\ran{f}$ is bounded by some $\hat{T}_{\mu_*}$-writable ordinal $\xi$. 
Now consider the following $\hat{T}_{\mu_*}$-computable function $g\colon \omega \rightarrow \mu_*$. Let $g(n)$ denote the halting position of the $n$th program, if this halts before time $\xi$, and $g(n)=0$ otherwise. 
Thus $L_{\hat{\lambda}_{\mu_*}}\models \cof(\mu_*)=\omega$. 
But Lemmas \ref{thetaomega equals theta<} and \ref{all values are equal} imply that $L_{\hat{\lambda}_{\mu_*}}\models \cof(\mu_*)>\omega$. 

To see that $\mu_*\leq \delta$, take any $\alpha<\mu_*$. 
Since $\mu_*=\mu_\omega$ by Lemma \ref{thetaomega equals theta<}, $T_\alpha$ can write an $\omega$-code for $\alpha$. 
Therefore, $T_\alpha$ can reach all cells by counting through this code. 
\end{proof} 




We call a cell \emph{eventually $T_\alpha$-reachable} if the head on the output tape eventually stabilizes on this cell. 
It is natural to ask whether a similar result holds for this notion of reachability. 
Let $\eta$ denote the least ordinal such that $T_\eta$ does not eventually reach all its cells. 

\begin{prop} 
$\delta=\eta$. 
\end{prop} 
\begin{proof} 
It is clear that $\delta\leq\eta$. 
Assume towards a contradiction that $\delta<\eta$. 
Then every cell of $T_\delta$ is eventually reachable. We partition the tapes into $\delta$ many portions of length $\delta$. For each cell $\xi$, we work in the $\xi$th portion and enumerate \emph{$\xi$-candidates} $(n,\alpha)$ that consist of a natural number and an ordinal by accidentally writing them via $\UU_\delta$. While the current $\xi$-candidate is considered, we pause $\UU_\delta$ and run the $n$th program on the $\xi$th portions of the tapes as long as the head position on the output tape is stable at the $\xi$th cell from time $\alpha$ onwards, with a code for $n$ on the output tape. Once the head moves, we run $\UU_\delta$ for the next step and switch to the next $\xi$-candidate. Note that if the $n$th program stabilizes at all, then it does so at or before time $\hat{\zeta}_\delta$ by Lemma \ref{loop}. This is accidentally writable by Lemma \ref{lambda zeta Sigma with and without parameters}. 
Thus the program eventually writes an output from which we can read off a function $f\colon \delta\rightarrow\omega$ mapping $\xi$ to $n$ as above. It is injective, since the $n$th program has a unique eventual head position, if its head stabilizes at all. 
Since $f\in L_{\hat{\zeta}_\delta}$, we have $L_{\hat{\zeta}_\alpha}\models |\delta|=\omega$. Since $L_{\hat{\lambda}_\alpha}\prec_{\Sigma_1} L_{\hat{\zeta}_\alpha}$, $L_{\hat{\lambda}_\alpha}\models |\delta|=\omega$.
But this contradicts Lemma \ref{thetaomega equals theta<} and Proposition \ref{delta equals mu}. 
\end{proof}

It is easy to see that the $T_\delta$-reachable cells form an interval, since $T_\delta$ can simulate $T_\alpha$ for all $T_\delta$-reachable $\alpha<\delta$ and $T_\alpha$ reaches all its cells. Hence $\delta$ equals the least $\alpha$ such that the $T_\alpha$-reachable cells are bounded. 


\begin{obs} 
There are arbitarily large countable ordinals $\alpha$ such that $T_\alpha$ can reach unboundedly many cells, but not all of them. 
\end{obs} 
\begin{proof} 
Recall that that the tape length is always assumed to be multiplicatively closed (see Section \ref{section setting}). 

We first claim that for any limit ordinal $\xi$ and any $i\in\omega$, the $\xi^{i}$th cell is $T_{\xi^{\omega}}$-reachable. (Note that $\xi^{\omega}$ is multiplicatively closed.) 
To see this, note that it is easy to implement a $T_\alpha$-program for ordinal multiplication (uniformly in $\alpha$) that sets the head to position $\beta\gamma$ when the cells with indices $\beta$ and $\gamma$ are marked. 
in this way, for any ordinal $\beta<\alpha$ and $i\in\omega$, we can move the head to position $\beta^{i}$ and write $1$s to the first $\beta^{i}$ many cells if the $\beta$th tape cell is marked with $1$ at the beginning of the computation.
By carrying out these procedures one after the other, we can also write $1$s to the first $\beta^{\omega}$ many cells and move the head $\beta^{\omega}$ many positions to the right. Now, to identify $\xi$ when $\alpha=\xi^{\omega}$, carry out this procedure successively for all $\nu<\alpha$, starting with $\nu=0$. As long as $\nu<\xi$, this will still leave $0$s on the tape, but when considering $\nu=\xi$, the whole tape will be filled with $1$s, which can be detected. Thus, we can identify the $\xi$th tape cell. Using the multiplication algorithm, it is now easy to see that $\xi^{i}$ is also reachable for all $i\in\omega$.

We further claim that that $\hat{\lambda}_{\xi}=\hat{\lambda}_{\xi^{\omega}}$ for all multiplicatively closed ordinals $\xi$. 
To see this, it suffices to note that tapes of length $\xi^{\omega}$ can be simulated on tapes of length $\xi$ by splitting the tape into $\omega$ many portions and simulating a tape of length $\xi^{i}$ on the $i$th portion. 

Now take any countable ordinal $\xi$ such that $L_{\hat{\lambda}_\xi+1}\models \xi$ is uncountable. 
There are unboundedly many such countable ordinals, since any image of $\omega_1$ in the transitive collapse of a countable elementary substructure of $L_{\omega_2}$ is of this form. 
Since we have seen that there are unboundedly many $T_{\xi^{\omega}}$-reachable cells, it remains to show that not all cells are $T_{\xi^{\omega}}$-reachable. 
Assuming otherwise, $L_{\hat{\lambda}_{\xi^\omega}+1}=L_{\hat{\lambda}_\xi+1}$ contains a surjection $f\colon \omega\rightarrow \xi^\omega$, contradicting the fact that  $\xi$ is uncountable in $L_{\hat{\lambda}_\xi+1}$. 
\end{proof}

\subsection{Writability strength} \label{writability strength} 

The next result shows that writability strength can decrease when the tape length increases. 

\begin{prop} \label{equivalent definitions of delta}
The following property of an ordinal $\alpha$ occurs for the first time at $\delta$: for some $\mu,\nu$ with $\omega\leq\mu\leq\nu<\alpha$, there is a $T_\nu$-writable but not $T_\alpha$-writable subset of $\mu$. 
\end{prop} 
\begin{proof} 
To see that $\delta$ has the required property, it suffices to find a $T_\nu$-writable but not $T_\delta$-writable subset of $\omega$ for some $\nu<\delta$. 
Note that $T_\nu$ can write an $\omega$-code of $\nu$ for all $\nu<\delta$ by Lemma \ref{thetaomega equals theta<} and Proposition \ref{delta equals mu}. 
Assuming the claim fails, $T_\delta$ could thus write an $\omega$-code for any $\nu<\delta$ and would therefore reach its $\nu$th cell. 

That $\delta$ is least follows from the fact that smaller devices can reach all their cells and therefore simulate all devices smaller than they are.
\end{proof} 


This suggests the question whether the writability strength for subsets of $\omega$ is comparable for different machines. 
The next result shows that this is the case. 


\begin{prop} \label{writability strength is comparable} 
For every $\alpha$, there is an ordinal $\tau_\alpha\leq\lambda_\alpha$ such that the $T_\alpha$-writable subsets of $\omega$ are exactly those contained in $L_{\tau_\alpha}$. Hence $T_\alpha$ and $T_\beta$ are comparable in their writability strength for subsets of $\omega$ for all $\alpha$, $\beta$. 
\end{prop}
\begin{proof} 
We first claim that every $T_{\alpha}$-writable real $x$ is contained in some $L_{\beta}$ with a $T_{\alpha}$-writable $\omega$-code. 
Note that $x\in L_{\lambda_{\alpha}}$ by Lemma \ref{sup of clockables equals sup of writables}. 
If $\beta$ is least with $x\in L_\beta$, then $L_{\beta}$ has a real code in $L_{\beta+1}$ by acceptability of the $L$-hierarchy. 
Hence  such a code is $T_{\alpha}$-accidentally writable without parameters. 
We run the universal $T_\alpha$-program $\UU_\alpha$ to search for an $\omega$-code of an $L$-level that contains $x$. Eventually, such an $\omega$-code for some $L_\tau$ is written on the output tape and the machine stops. 

It remains to see that every real in some $L_\tau$ with a $T_\alpha$-writable $\omega$-code $y$ for $L_\tau$ is itself $T_\alpha$-writable, but this is clear since each element of $L_\tau$ is coded in $y$ by a natural number. 
\end{proof}

It is easy to see that the previous result fails for subsets of other ordinals if the machine has non-reachable cells.

We now turn to characterizations of $\delta$ via eventually and accidentally writable sets. The next result follows from Proposition \ref{delta equals mu}, the fact that every $T_\alpha$-accidentally writable subset of $\alpha$ is an element of $L_{\Sigma_\alpha}$ by Lemma \ref{loop} and the discussion at the end of Section \ref{subsection L-levels}. 

\begin{prop} 
The following properties of $\alpha$ occur first at $\delta$: 
\begin{enumerate}[(a)] 
\item \label{condition for eventual writability 1} 
There is no $T_\alpha$-eventually writable $\xi$-code ($\omega$-code) of $\alpha$ for some $\xi<\alpha$. 
\item \label{condition for eventual writability 1} 
There is no $T_\alpha$-accidentally writable $\xi$-code ($\omega$-code) of $\alpha$ for some $\xi<\alpha$. 
\end{enumerate} 
\end{prop} 

\noindent
We say that a set of ordinals \emph{has a gap} if it is not an interval. 
For standard \ITTM s there are no gaps in the writable ordinals, since from a code for an ordinal one can write a code for any smaller ordinal by simply truncating the code \cite[Theorem 3.7]{HL}. However, for $\delta$-codes truncating would require addressing every tape cell, which is not possible when there are non-reachable cells. 

\begin{lem} 
$\delta$ is least such that the $T_\delta$-writable ordinals have a gap. 
\end{lem} 
\begin{proof} 
There are no gaps in the $T_\alpha$-writable ordinals for $\alpha<\delta$, since every cell is reachable and hence codes can be truncated at any length. 
We now show that $[\theta,\delta)$ is the first gap for $T_\delta$, where $\theta$ is the least cell that is not $T_\delta$-reachable. To see this, note that it follows from the equality $\delta=\mu_*=\mu_\omega$ in Lemma \ref{thetaomega equals theta<} and Proposition \ref{delta equals mu} that every $T_\delta$-reachable $\alpha$ has a $T_\delta$-writable $\omega$-code and it is also clear that $\delta$ has a $T_\delta$-writable $\delta$-code. If some $\alpha\in [\theta,\delta)$ had a $T_\delta$-writable $\delta$-code, then one would be able to reach $\alpha$ by counting through the code, but this contradicts the choice of $\theta$. 
\end{proof}

\subsection{The role of parameters} \label{subsection role of parameters} 

While $\Sigma_\alpha=\hat{\Sigma}_\alpha$ by Lemma \ref{lambda zeta Sigma with and without parameters}, the next result shows that the analogous statement for $\hat{\lambda}_\alpha$ and $\hat{\zeta}_\alpha$ fails.  
We would like to thank Philip Welch for providing a proof of the implication from \ref{role of parameters 1} to \ref{role of parameters 3}. 
This answered an open question in a preliminary version of this paper. 

\begin{thm} \label{role of parameters} 
The following statements are equivalent: 
\begin{enumerate}[(a)]
\item \label{role of parameters 1} 
$\alpha$ is countable in $L_{\lambda_\alpha}$. 
\item \label{role of parameters 2} 
As in \ref{role of parameters 1}, but with $\lambda_\alpha$ replaced by $\hat{\lambda}_\alpha$, $\zeta_\alpha$, $\hat{\zeta}_\alpha$ or $\Sigma_\alpha=\hat{\Sigma}_\alpha$. 
\item \label{role of parameters 3}
$\hat{\lambda}_\alpha=\lambda_\alpha$. 
\item \label{role of parameters 4}
$\hat{\zeta}_\alpha=\zeta_\alpha$. 
\end{enumerate} 
\end{thm} 
\begin{proof} 
To see the equivalence of \ref{role of parameters 1} and  \ref{role of parameters 2}, recall that $L_{\hat{\lambda}_\alpha}\prec_{\Sigma_1} L_{\hat{\zeta}_\alpha}\prec_{\Sigma_1} L_{\hat{\Sigma}_\alpha}=L_{\Sigma_\alpha}$ by Theorem \ref{lambda-zeta-Sigma submodel characterization} and $L_{\lambda_\alpha}\prec_{\Sigma_1^{(\alpha)}} L_{\zeta\alpha}\prec_{\Sigma_1^{(\alpha)}} L_{\Sigma_\alpha}$, where $\Sigma_1^{(\alpha)}$ denotes $\Sigma_1$-formulas only in the parameter $\alpha$, as discussed at the end of Section \ref{section writable and clockable ordinals}. 
Thus the claim follows from the fact that countability of $\alpha$ is expressible by a $\Sigma_1^{(\alpha)}$-formula. 

Now assume \ref{role of parameters 1} and  \ref{role of parameters 2}. 
Since the set of ordinals with $T_\alpha$-writable $\omega$-codes is downwards closed, $\alpha$ has a $T_\alpha$-writable $\omega$-code. 
Then any cell is $T_\alpha$-reachable by counting through the code. Hence \ref{role of parameters 3} and \ref{role of parameters 4} hold. 

Conversely, assume that \ref{role of parameters 1} and  \ref{role of parameters 2} fail. 
Thus $\alpha$ is uncountable in $L_{\Sigma_\alpha}$. 

To show that \ref{role of parameters 3} fails, let $H$ denote the set of $e\in \NN$ such that $T_\alpha[e]$ halts and outputs an $\alpha$-code for an ordinal $\gamma_e$. 
Since all $T_\alpha$-clockable ordinals are below $\lambda_\alpha$ by Lemma \ref{sup of clockables equals sup of writables}, $H\in L_{\lambda_\alpha+1}$. 
Since $\alpha$ is uncountable in $L_{\Sigma_\alpha}$ by our assumption and the $L$-hierarchy is acceptable \cite[Theorem 1]{MR0239977}, $H\in L_\alpha\subseteq L_{\lambda_\alpha}$. 
Note that the function $f\colon H\rightarrow \lambda_\alpha$, $f(e)=\gamma_e$ is $\Sigma_1$-definable from $H$ over $L_{\hat{\lambda}_\alpha}$. 
Since $\hat{\lambda}_\alpha$ is admissible by Lemma \ref{lambdaalpha is admissible}, $\lambda_\alpha=\sup_{e\in H}\gamma_e<\hat{\lambda}_\alpha$. 

To show that \ref{role of parameters 4} fails, let $H^*$ denote the set of $e\in\NN$ such that $T_\alpha[e]$ eventually outputs an $\alpha$-code for an ordinal $\gamma^*_e$. 
Since all stabilization times are below $\zeta_\alpha$ by Lemma \ref{sup of stabilization times equals sup of writables with parameters}, $H^*\in L_{\zeta_\alpha+1}$ and therefore $H^*\in L_\alpha \subseteq L_{\zeta_\alpha}$, as above. 
Moreover, the function $g\colon H^*\rightarrow \zeta_\alpha$, $g(e)=\gamma^*_e$ is $\Sigma_2$-definable from $H^*$ over $L_{\hat{\zeta}_\alpha}$. 
Since $\hat{\zeta}_\alpha$ is $\Sigma_2$-regular by Lemma \ref{lambdaalpha is admissible}, $\zeta_\alpha=\sup_{e\in H^*}\gamma^*_e<\hat{\zeta}_\alpha$. 
\end{proof} 

\noindent
We obtain the next Corollary via Lemma \ref{equivalent definitions of mugamma} and Proposition \ref{delta equals mu}.  

\begin{cor} 
$\delta$ equals the least ordinal $\alpha$ with each of the following properties: 
\begin{enumerate}[(a)] 
\item 
$\lambda_\alpha<\hat{\lambda}_\alpha$. 
\item 
$\zeta_\alpha<\hat{\zeta}_\alpha$. 
\end{enumerate} 
\end{cor}

\subsection{Upper and lower bounds} \label{subsection upper and lower bounds} 


We have the following upper bound for $\delta$. 
Let $\sigma$ be the least ordinal $\alpha$ such that every $\Sigma_1$-statement true in $L$ already holds in $L_\alpha$ ($L_\sigma$ equals the $\Sigma_1$-hull of $\emptyset$ in $L$, since it contains every uniquely $\Sigma_1$-definable set). 
Since both the statement that a program halts and the existence of an ever-repeating loop are $\Sigma_1$-statements, the existence of $\delta$ is a $\Sigma_1$-statement and hence $\delta<\sigma$. 



For a lower bound, we see that $\delta$ is a closure point of the function mapping $\alpha$ to $\Sigma_\alpha$. 

\begin{thm} \label{Sigmaalpha is less than delta}
$\Sigma_\alpha<\delta$ for all $\alpha<\delta$.\footnote{This strengthens the result from \cite{MR3210081} that $\zeta<\delta$ and an unpublished result by Robert Lubarsky that $\Sigma<\delta$.} 
\end{thm} 
\begin{proof} 
By Proposition \ref{delta equals mu}, Lemma \ref{thetaomega equals theta<}, Lemma \ref{all values are equal} and the discussion after it, $\delta$ is a regular cardinal in the admissible set $L_{\hat{\lambda}_\delta}$. Hence there is a strictly increasing sequence $\langle \xi_\beta\mid\beta<\delta\rangle\in L_{\hat{\lambda}_\delta}$ of ordinals with $\alpha<\xi_\beta<\delta$ such that $\langle L_{\xi_\beta}\mid \beta<\delta\rangle\in L_{\hat{\lambda}_\delta}$ is a chain of elementary substructures of $L_\delta$. 
In particular, $L_{\xi_0}\prec_{\Sigma_1} L_{\xi_1}\prec_{\Sigma_2} L_{\xi_2}$. Since the triple $(\hat{\lambda}_\alpha, \hat{\zeta}_\alpha, \hat{\Sigma}_\alpha)$ is least with this property by Theorem \ref{lambda-zeta-Sigma submodel characterization}, we have $\Sigma_\alpha=\hat{\Sigma}_\alpha\leq \xi_2<\delta$. 
\end{proof}

\section{Open questions}


Since we considered various conditions that occur at $\delta$ for the first time, it is natural to ask which of them are equivalent everywhere. 

\begin{qu} 
Which of the conditions in Theorem \ref{main theorem} are equivalent for all ordinals? 
\end{qu} 

Throughout the paper, we worked with the functions mapping a multiplicatively closed ordinal $\alpha$ to the values $\lambda_\alpha$, $\hat{\lambda}_\alpha$, $\zeta_\alpha$, $\hat{\zeta}_\alpha$ and $\Sigma_\alpha=\hat{\Sigma}_\alpha$. 
The versions of these functions with parameters are monotone, since $\hat{T}_\beta$ can simulate $\hat{T}_\alpha$ for $\alpha\leq\beta$. 
Moreover, the versions without parameters are monotone below $\delta$ for a similar reason, and at $\delta$ by Theorem \ref{Sigmaalpha is less than delta}. 

\begin{qu} 
Are the functions $\alpha\mapsto \lambda_\alpha, \zeta_\alpha$ monotone above $\delta$? 
\end{qu} 





We are further interested in the supremum $\theta_\alpha$ of $T_\alpha$-reachable cells. 
For instance, one can ask the next question. 

\begin{qu}
Is the function $\alpha\mapsto \theta_\alpha$ monotone? 
\end{qu}





Finally, we ask whether similar results to those in this paper hold for 
machines with $\Sigma_n$-limit rules \cite{FW}.

\bibliographystyle{alpha}
\bibliography{references}

\begin{thebibliography}{COW18}

\bibitem[BP68]{MR0239977}
George Boolos and Hilary Putnam.
\newblock Degrees of unsolvability of constructible sets of integers.
\newblock {\em J. Symbolic Logic}, 33:497--513, 1968.

\bibitem[COW18]{CO}
Merlin Carl, Sabrina Ouazzani, and Philip~D. Welch.
\newblock Taming {K}oepke's {Z}oo.
\newblock In {\em CiE}, 2018.

\bibitem[FW11]{FW}
Sy-David Friedman and Philip Welch.
\newblock Hypermachines.
\newblock {\em J. Symbolic Logic}, 76(2):620--636, 2011.

\bibitem[HL00]{HL}
Joel~David Hamkins and Andy Lewis.
\newblock Infinite time {T}uring machines.
\newblock {\em J. Symbolic Logic}, 65(2):567--604, 2000.

\bibitem[Koe05]{K1}
Peter Koepke.
\newblock Turing computations on ordinals.
\newblock {\em Bull. Symbolic Logic}, 11(3):377--397, 2005.

\bibitem[Koe09]{K}
Peter Koepke.
\newblock Ordinal computability.
\newblock In {\em Mathematical theory and computational practice}, volume 5635
  of {\em Lecture Notes in Comput. Sci.}, pages 280--289. Springer, Berlin,
  2009.

\bibitem[KS09]{KS}
Peter Koepke and Benjamin Seyfferth.
\newblock Ordinal machines and admissible recursion theory.
\newblock {\em Ann. Pure Appl. Logic}, 160:310--318, 2009.

\bibitem[Rin14]{MR3210081}
Benjamin Rin.
\newblock The computational strengths of {$\alpha$}-tape infinite time {T}uring
  machines.
\newblock {\em Ann. Pure Appl. Logic}, 165(9):1501--1511, 2014.

\bibitem[Wel00]{MR1734198}
P.~D. Welch.
\newblock The length of infinite time {T}uring machine computations.
\newblock {\em Bull. London Math. Soc.}, 32(2):129--136, 2000.

\bibitem[Wel09]{MR2493990}
P.~D. Welch.
\newblock Characteristics of discrete transfinite time {T}uring machine models:
  halting times, stabilization times, and normal form theorems.
\newblock {\em Theoret. Comput. Sci.}, 410(4-5):426--442, 2009.

\end{thebibliography}

\end{document}